\documentclass[12pt,a4paper]{amsart}

\usepackage{amsmath,amsthm,amssymb,latexsym,a4wide,color}

\newtheorem{theorem}{Theorem} [section]

\newtheorem{corollary}[theorem]{Corollary}

\theoremstyle{definition}

\numberwithin{equation}{section}

\newcommand{\mcL}{\mathrel{\mathcal{L}}}
\newcommand{\mcR}{\mathrel{\mathcal{R}}}

\newcommand{\tmcLU}{\mathrel{\widetilde{\mcLx}_U}}
\newcommand{\tmcLV}{\mathrel{\widetilde{\mcLx}_V}}

\newcommand{\mcLx}{\mathcal{L}}
\newcommand{\mcRx}{\mathcal{R}}
\newcommand{\mcHx}{\mathcal{H}}
\newcommand{\mcDx}{\mathcal{D}}
\newcommand{\tmcLUx}{\widetilde{\mcLx}_U}
\newcommand{\tmcRUx}{\widetilde{\mcRx}_U}

\newcommand{\mcS}{\mathcal{S}}
\newcommand{\mcG}{\mathcal{G}}
\newcommand{\mcC}{\mathcal{C}}
\newcommand{\mcO}{\mathcal{O}}
\newcommand{\mcU}{\mathcal{U}}
\newcommand{\mcK}{\mathcal{K}}

\newcommand{\dd}{\mathbf{d}}
\newcommand{\rr}{\mathbf{r}}

\newcommand{\oml}{\mathrel{\omega^{\ell}}}
\newcommand{\omr}{\mathrel{\omega^{r}}}
\newcommand{\om}{\mathrel{\omega}}
\newcommand{\omlx}{\omega^{\ell}}
\newcommand{\omrx}{\omega^{r}}
\newcommand{\omx}{\omega}
\newcommand{\eps}{\varepsilon}

\newcommand{\muc}{\mu^\circ}

\newcommand{\IG}{\operatorname{IG}}
\newcommand{\RIG}{\operatorname{RIG}}
\newcommand{\GL}{\operatorname{GL}}
\newcommand{\Ker}{\operatorname{Ker}}

\newcommand{\sthe}{(S,\theta)}

\def\act#1#2{{^{#1}\kern -2pt {#2}}}
\def\acth#1#2{{^{#1}\kern -1pt {#2}}}

\def\softd{{\leavevmode\setbox1=\hbox{d}%
     \hbox to 1.05\wd1{d\kern-0.4ex{\char039}\hss}}}

\title[Regular semigroups]{Structure theory of regular semigroups}

\author{M\'aria B.\ Szendrei}
\address{Bolyai Institute, University of Szeged, Aradi v\'ertan\'uk tere 1, H-6720 Szeged, Hungary; Alfr\'ed R\'enyi 
Institute of Mathematics, Hungarian Academy of Sciences, Re\'altanoda utca 13--15, H-1053 Budapest, Hungary}
\email{m.szendrei@math.u-szeged.hu}
\thanks{The author was partially supported by the National Research, Development and Innovation Office, grants K115518 and K128042.}

\begin{document}

\date{16 August, 2019}

\begin{abstract} 
This survey aims to give an overview of several substantial developments of the last 50 years in the structure theory of regular semigroups and to shed light on their impact on other parts of semigroup theory.
\end{abstract}

\maketitle

\section{Introduction}

Research interest has been centred around regular semigroups from the very beginning when semigroups appeared as independent algebraic structures to study.
In the earliest major result on semigroups, Su\v{s}kevi\v{c} (1928) described the structure of finite simple semigroups which form an important class of regular semigroups.
After the individual viewpoints and methods of the theory of regular semigroups had been developed, and the structure of various special classes of regular semigroups had been described, the time came in the late 1960's and in the 1970's to focus on the general structure of regular semigroups. 
Inverse semigroups form a prominent subclass of regular semigroups which have applications in a number of areas of mathematics, and also outside mathematics: differential geometry, theory of $C^*$-algebras, combinatorial group theory, model theory, linear logic, tilings, quasicrystals and solid-state physics (\cite{Lawson}).
The significance of inverse semigroups in the structure theory of regular semigroups is due to the fact that the structure of inverse semigroups is much simpler than that of regular semigroups in general, and so the results proved for inverse semigroups serve as initial steps for proving (more) general results for regular semigroups.

In the last 50 years, a huge number of papers have been published on the structure of regular semigroups. 
A high diversity of subclasses 
(e.g., completely regular semigroups, orthodox semigroups, locally inverse semigroups, $E$-solid (also called quasi-orthodox) semigroups, regular $*$-semigroups, $P$-regular semigroups, regular semigroups with inverse trans\-versals) 
have been investigated, and many new ideas and methods have been found.
In order to keep the length of the paper in a reasonable range, 
we had to make some choices what to cover in this survey.
For example, we will not cover completely regular semigroups, because
due to their special features, their structure theory has its own methods and constructions.
For the results on the structure of completely regular semigroups, we refer the reader to the monograph by Petrich and Reilly \cite{Petrich-Reilly}.
We will omit structure theorems where specific constructions are introduced to describe the structure of certain regular semigroups.
The reason for this choice is partly that these constructions are often long and technical, and partly that
a large portion of these constructions are special cases of extensions of regular semigroups
introduced and studied in detail by Pastijn and Petrich \cite{Pastijn-Petrich-85}. 
Instead of such structure theorems, we present in the paper embedding theorems where regular semigroups are embedded into products which generalize semidirect products of groups.
Finally, we will not discuss the theory of existence varieties despite the fact that bifree objects and bi-identities of existence varieties have found applications in the structure theory of regular semigroups. 
For an introduction to the theory of existence varieties, the reader is referred to the survey papers by 
Auinger \cite{Auinger-X02}, Jones \cite{Jones-95} and Trotter \cite{Trotter-96}. 

The topics of the structure theory of regular semigroups that we will cover in this survey were inspired by the three major approaches that describe the structure of all inverse semigroups.
One of them is Munn's approach via fundamental inverse semigroups (1966-1970), and another one is due to McAlister and Lawson via semidirect products of semilattices by groups (1974-76, 1992).
Both build up inverse semigroups, in a certain sense, from semilattices and groups.
The third approach found by Schein constructs inverse semigroups from ordered groupoids over semilattices and vice versa (1965), but the idea of applying categories goes back to Ehresmann (1958).
The last two approaches form the basis for the algebraic tools that are applied in Lawson's monograph \cite{Lawson} to present the interrelationships among inverse semigroups, partial symmetries and global symmetries.

The three approaches mentioned in the previous paragraph inspired intensive research on regular semigroups, and several outstanding results were achieved in the 1970's that have had far reaching influence in the theory of regular semigroups and beyond.
Section \ref{sect:fund} contains generalizations of Munn's approach to regular semigroups, and gives an account of an extension of Nambooripad's approach via biordered sets to arbitrary semigroups, and of the topic of idempotent generated semigroups.
Section \ref{sect:grd-crc} is mainly about Nambooripad's generalization of Schein's result for the class of regular semigroups, but his other approach via cross-connections and several recent results are also presented.
In contrast to Munn's approach and Schein's result, the McAlister--Lawson approach has not been extended to the whole class of regular semigroups.
The results obtained for orthodox and for locally inverse semigroups form the main topic of Section \ref{sect:McA-Law}.
They suggest that a generalization for the whole class of regular semigroups, if it exists at all, will require dissimilar approaches and methods.
The section also contains works motivated, at least partly, by the study of the main topic.

\section{Preliminaries} \label{Prelim}

The reader is assumed to be familiar with the usual terminology and notation and with the basic facts of the theory of regular semigroups, or he is referred for background to the monographs by Lawson \cite{Lawson}, Howie \cite{Howie-76,Howie-95}, Grillet \cite{Grillet-95}, and Petrich \cite{Petrich}.
In several cases, terminology and notation differ in them.
We prefer the commonly used expressions `fundamental inverse semigroup', `locally inverse semigroup' and `regular idempotent generated semigroup' to `antigroup' (\cite{Petrich}), `pseudoinverse semigroup' (\cite{Grillet-95}) and `semiband' (\cite{Howie-95}), respectively.
Unlike in \cite{Lawson} and partly in \cite{Grillet-95}, full and partial mappings (functions) are written as right operators, and their composition as well as composition of morphisms in categories are formed from left to right.
The set of idempotents of a semigroup $S$ is denoted $E(S)$, and the maximum idempotent separating congruence on a regular semigroup by $\mu$.

We remark that, although the notion of a $\lambda$-semidirect product of $K$ by $T$ is defined in \cite{Lawson} only for the case where both $K$ and $T$ are inverse semigroups, the construction works for any semigroup $K$, and if $K$ is regular then so is the $\lambda$-semidirect product of $K$ by $T$.

Let $S$ be a regular semigroup and $\theta$ a congruence on $S$.
Recall that, in general, a $\theta$-class containing an idempotent element is a subsemigroup of $S$ that need not be regular.
If each $\theta$-class containing an idempotent element belongs to a given class $\mcK$ of regular semigroups then we say 
$\theta$ is a {\em congruence over $\mcK$}.
In particular, if $\theta$ is an inverse semigroup congruence (more specially, a group congruence) then each $\theta$-class containing an idempotent element is regular, and the union of these $\theta$-classes is also a regular subsemigroup of $S$.
The latter subsemigroup is called the {\em kernel of $\theta$} and is denoted $\Ker\theta$.
For example, a congruence is idempotent separating if and only if it is over groups, and it is idempotent pure if and only if it is over bands.

Let $T$ be a regular semigroup and $\mcK$ a class of regular semigroups.
If $S$ is a regular semigroup and $\theta$ is a congruence on $S$ (over $\mcK$) such that $S/\theta$ is 
isomorphic to $T$ then $\sthe$ is called an {\em extension by $T$ (over $\mcK$)}, or a coextension of $T$ (over $\mcK$).
In this survey, we use the former terminology.
Now let $K$ be a regular semigroup and $T$ an inverse semigroup (in particular, a group). 
If $S$ is a regular semigroup and $\theta$ is a congruence on $S$ such that $S/\theta$ is isomorphic to $T$ and $\Ker\theta$ is isomorphic to $K$ then $\sthe$ is called an {\em extension of $K$ by $T$}. 
Note that if $S$ is an inverse semigroup then $\sthe$ is an extension of $K$ by $T$ if and only if $(\iota,S,\theta^\natural)$, where $\iota$ is the inclusion mapping $\Ker\theta\to S$, is a normal extension of $K$ by $T$ along 
$\iota\theta^\natural$ (\cite[Section 5.1]{Lawson}, \cite[Section VI.6]{Petrich}).
If $\sthe$ and $(S',\theta')$ are extensions then an injective homomorphism $\phi\colon S\to S'$ is called an {\em embedding of $\sthe$ into $(S',\theta')$} if $\theta=\ker(\phi\theta'^\natural)$.
If a construction (e.g., semidirect product, $\lambda$-semidirect product) which produces a regular semigroup 
$K\star T$ from regular semigroups $K$ and $T$ on a subset of the set $K\times T$ such that the second projection 
$\pi\colon K\star T\to T$ is a surjective homomorphism then $(K\star T,\ker\pi)$ is an extension by $T$, and if $T$ is an inverse semigroup then $(K\star T,\ker\pi)$ is an extension of $\Ker(\ker\pi)$ by $T$. 
In the context of extensions, we understand this extension when referring simply to $K\star T$.

By the {\em core $C(S)$} of a semigroup $S$ we mean the subsemigroup of $S$ generated by $E(S)$ provided $E(S)$ is not empty.
It is worth mentioning that if $S$ is regular then so is $C(S)$ (\cite[Section II.4, Exercise 15]{Howie-76}, \cite[Exercise 2.6.15]{Howie-95}).
The {\em self-conjugate core $C_\infty(S)$} of a regular semigroup $S$ is defined to be the minimum subsemigroup in $S$ containing $C(S)$ and being closed under conjugation. 
Here a conjugate of $a\in S$ is an element $x'ax$ where $x\in S^1$ and $x'$ is an inverse of $x$. 
Notice that $C_\infty(S)$ is contained in the kernel of the least inverse semigroup congruence on $S$. 
An {\em $E$-solid semigroup $S$} is defined to be a regular semigroup such that $C(S)$ is completely regular.
It was shown by Hall \cite[Supplement]{Yamada-79} and Trotter \cite{Trotter-91} 
that if $S$ is $E$-solid then $C_\infty(S)$ is also completely regular, and  $C_\infty(S)$ is equal to the kernel of the least inverse semigroup congruence on $S$.

\section{Fundamental regular semigroups} \label{sect:fund}

In the 1970's, the general structure of regular semigroups was a central topic of research, and several different approaches appeared.
In this section we outline the major results on the structure of regular semigroups
(Sections \ref{Hall}--\ref{Namboo-fund}), mostly achieved in the 1970's, which generalize Munn's approach to the description of inverse semigroups via fundamental inverse semigroups (\cite[Theorems 5.2.7--5.2.9]{Lawson}).
Due to Easdown \cite{Easdown}, Nambooripad's concept of a biordered set has turned out to play significant role among arbitrary semigroups.
Nambooripad's results on regular idempotent generated semigroups (\cite{Nambooripad-79}), combined with this result, 
have been the inspiration for
the investigation of the structure of free idempotent generated semigroups.
Several significant results of this topic are mentioned in Section \ref{id-gen} together with ones on the structure of the biordered sets of members of important classes of semigroups.
At the end of the section we sketch several results which extend those on fundamental regular semigroups to wider classes of semigroups.

In an inverse semigroup $S$, the set $E(S)$ of the idempotent elements forms a subsemilattice --- a nice structure which is easy to work with ---, and among regular semigroups, this property characterizes them.
It was observed by Munn 
that  if $E$ is a semilattice then the set $T_E$ of all isomorphisms between the 
principal ideals of $E$ forms an inverse subsemigroup of the symmetric inverse semigroup on $E$. 
The inverse semigroup $T_E$ and all of its full inverse subsemigroups are fundamental.
Moreover, conjugation by an element of an inverse semigroup $S$ determines an element of $T_{E(S)}$,
thus defining a representation of $S$ in $T_{E(S)}$ whose kernel is the maximum idempotent separating congruence on $S$.
Consequently, the fundamental inverse semigroups with semilattice $E$ are, up to isomorphism, just the full inverse subsemigroups of $T_E$, and any inverse semigroup with semilattice $E$ is an idempotent separating extension by a full inverse subsemigroup of $T_E$.
In other words, the latter statement says that any inverse semigroup with semilattice $E$ is an extension by a full inverse subsemigroup of $T_E$ over groups.
Ultimately, any inverse semigroup can be built up from semilattices and groups.
The inverse semigroup $T_E$ and the representation mentioned are called the 
\emph{Munn semigroup of $E$} and the \emph{Munn representation of $S$}, respectively.

The structure of a regular semigroup is much more complicated than that of an inverse semigroup:
if $S$ is a regular semigroup then its elements might have several inverses, the natural partial order is not necessarily compatible with multiplication, the set $E(S)$ of idempotents might not be closed under multiplication, etc. 
Therefore it is by no means clear how to generalize the Munn representation --- or more generally, a theorem on the structure of inverse semigroups --- for regular semigroups.
It turns out, however, that, similarly to the inverse case, the behavour of a regular semigroup is strongly influenced by the structure of its idempotents.

\subsection{Hall's approach: fundamental orthodox and regular semigroups via bands and idempotent generated semigroups} \label{Hall}

The first attempt to generalize the Munn representation, due to Hall \cite{Hall-70} (see also \cite[Section VI.2]{Howie-76}), 
focuses on the intermediate class of orthodox semigroups, that is, on regular semigroups $S$ such that $E(S)$ is a 
subband of $S$.
For any band $B$, he constructs a generalization $W_B$ of the Munn semigroup, called the Hall semigroup of $B$ (\cite{Howie-76}), with the property that the fundamental orthodox semigroups with band $B$ are just the full orthodox subsemigroups of $W_B$.
He also generalizes the Munn representation by defining, for any orthodox semigroup $S$, a homomorphism from $S$ into $W_{E(S)}$ such that its range is a full orthodox subsemigroup.
Similarly to the inverse case, this implies that each orthodox semigroup with band $B$ is an extensions by a full orthodox subsemigroup of $W_B$ over groups.
In a subsequent paper, Hall \cite{Hall-73} extends this result to regular semigroups where the role of a band is taken over by
a regular idempotent generated semigroup, and in particular, the role of $E(S)$ is taken over by the core $C(S)$ of $S$.
The regular semigroup $W_B$ constructed from a regular idempotent generated semigroup $B$ is fundamental, but its core is, in general, isomorphic to $B/\mu$ rather than to $B$ itself.
(For example, a fundamental completely simple semigroup is necessarily a rectangular band, but there exist non-fundamental idempotent generated completely simple semigroups.)
Consequently, the representation of a regular semigroup $S$ in $W_{C(S)}$ is also weaker than in the orthodox case.
In particular, if $B$ is fundamental then the fundamental regular semigroups with core $B$ are precisely the full regular subsemigroups of $W_B$. 
Note that both in the orthodox and in the regular cases, $W_B$ consists of pairs of transformations where the components of a pair are defined on the sets $B/\mcLx$ and $B/\mcRx$, respectively.

\subsection{Grillet's approach: fundamental regular semigroups via cross-connections} \label{Grillet}

Grillet \cite{Grillet-74a, Grillet-74b} (see also \cite[Sections VIII.1--2]{Grillet-95})
generalizes the Munn representation in such a way that the partially ordered sets $S/\mcLx$ and $S/\mcRx$ are used in the place of the semilattice of $S$.
Note that these partially ordered sets are order isomorphic to $B/\mcLx$ and $B/\mcRx$ applied in Hall's approach.
However, while Hall's construction of $W_B$ applies the multiplication of $B$, Grillet characterizes
the partially ordered sets $S/\mcLx$ and $S/\mcRx$ of regular semigroups $S$ in an abstract way, and call them 
\emph{regular partially ordered sets}. 
Furthermore, by introducing the notion of a \emph{cross-connection} between two regular partially ordered sets, he shows that 
a pair of regular partially ordered sets stems from a regular semigroup $S$ as $S/\mcLx$ and $S/\mcRx$ if and only if there is a cross-connection between them. 
This allows him to introduce, for any pair $I,\Lambda$ of regular partially ordered sets and for any cross-connection
$(\Gamma,\Delta)$ between them, a regular semigroup $T_{I,\Lambda;\Gamma,\Delta}$ such that this semigroup plays the same role among the regular semigroups $S$ with $S/\mcRx=I$, $S/\mcLx=\Lambda$ and induced cross-connection $(\Gamma,\Delta)$
that is played by the Munn semigroup $T_E$ among the inverse semigroups with semilattice $E$.
In particular, Munn's results can be deduced from Grillet's approach. 
Note also that \cite{Grillet-74c} contains further details on fundamental regular semigroups, and \cite{Grillet-74d} provides a general structure theorem for all regular semigroups, based on the fact that each regular semigroup is an extension by a fundamental one over groups.

\subsection{Nambooripad's approach: fundamental regular semigroups via biordered sets} \label{Namboo-fund}

The approach chosen by Nambooripad in \cite{Nambooripad-73} (see also \cite{Nambooripad-75a})
axiomatizes the structure of the set of idempotents of a regular semigroup as a regular biordered set, 
and generalizes the Munn representation for fundamental regular semigroups with a given biordered set. 
Note that this definition of a regular biordered set uses partial transformations. 
Instead of a regular biordered set, Clifford \cite{Clifford-75a} introduces the notion of a {\em warp}, which is a set equipped with a partial operation such that certain axioms are satisfied.
He defines the warp of a regular semigroup $S$ to be the set $E(S)$ with the partial operation induced by the operation of $S$, and generalizes the Munn representation in this framework.
It was Clifford's intention with his approach to generalize not only the Munn representation but also the Hall representation of orthodox semigroups \cite{Hall-70} (see Section \ref{Hall}).
However, in the general regular case, the same difficulties appear with this approach that were mentioned in Section \ref{Hall}: the warps of $S$ and of $S/\mu$ need not be isomorphic for a regular semigroup $S$.

Influenced by Clifford's idea, Nambooripad redefines a biordered set in \cite{Nambooripad-79} as a specific warp.
We summarize Nambooripad's results in this setting (see also \cite[Sections VIII.3--4]{Grillet-95}). 

If $S$ is any semigroup then 
let the quasi-orders $\omlx$ and $\omrx$ 
be defined, for every $e,f\in E(S)$, by
\[e\oml f\quad \hbox{if and only if}\quad ef=e,\quad \hbox{and}\quad e\omr f\quad \hbox{if and only if}\quad fe=e,\]
and consider the partial operation on $E(S)$ which is induced on the set
$\omlx\cup(\omlx)^{-1}\cup\omrx\cup(\omrx)^{-1}$ by the multiplication of $S$.
Notice that $e\oml f$ is equivalent to the inclusion $S^1e\subseteq S^1f$, and similarly, $e\omr f$ to $eS^1\subseteq fS^1$. 
Consequently, we have $\mcLx=\omlx\cap(\omlx)^{-1}$ and $\mcRx=\omrx\cap(\omrx)^{-1}$ in $E(S)$, and the relation 
$\omx=\omlx\cap\omrx$ is the natural partial order on $E(S)$.
The {\em sandwich set} of $e,f\in E(S)$ is defined as follows:
\begin{eqnarray*}
\lefteqn{\mcS(e,f)=\{h\in E(S): h\oml e,\ h\omr f,\ \hbox{and}\ }\\
&& \qquad\quad \hbox{if}\ g\in E(S)\ \hbox{with}\ g\oml e,\ g\omr f\ \hbox{then}\ eg\omr eh\ \hbox{and}\ gf\oml hf\}.
\end{eqnarray*}
The abstract notion of a \emph{biordered set} is defined to be a set equipped with a partial operation that
satisfies certain axioms.
These axioms are valid in the biordered set of every semigroup.
A biordered set is called \emph{regular} if $\mcS(e,f)$ is non-empty for any elements $e,f$.
In particular, semilattices are just the regular biordered sets where $\omlx=\omrx(=\omx)$, in which case each sandwich set is
necessarily a singleton, and the element of $\mcS(e,f)$ is the greatest lower bound of $e$ and $f$. 

Nambooripad constructs a fundamental regular semigroup $T_E$ which 
is an exact generalization of the Munn semigroup of a semilattice, and
generalizes the Munn representation in the most elegant form possible.

\begin{theorem}[Nambooripad \cite{Nambooripad-79}]\label{Namb-fund}
\begin{enumerate}
\item \label{Namb-fund1}
If $E$ is a regular biordered set then $T_E$ is a fundamental regular semigroup whose biordered set is $E$, 
and each full regular subsemigroup of $T_E$ is a fundamental regular semigroup.
\item
For every regular semigroup $S$ with biordered set $E$, there exists a representation 
$S\to T_E$ which preserves $E$ and whose kernel is the maximum idempotent separating congruence on $S$.
\end{enumerate}
\end{theorem}

In particular, statement (\ref{Namb-fund1}) justifies the notion of an abstract regular biordered set.

\begin{corollary}
Each regular biordered set is a biordered set of a regular semigroup.
\end{corollary}

\subsection{Biordered sets and idempotent generated semigroups} \label{id-gen}

There are many natural examples of idempotent generated semigroups, e.g., for every natural number $n$, the semigroup $T_n\setminus S_n$ of all singular transformations on the set $\{1,\ldots,n\}$ and the semigroup $M_n(F)\setminus \GL_n(F)$ of all singular $n\times n$ matrices over a field $F$. 
A consequence of this result on transformations is that each finite semigroup is embeddable in a finite idempotent generated semigroup, and the analogous statement holds in general (\cite[Corollary 6.3.3 and Theorem 6.3.4]{Howie-95}).
Moreover, in each semigroup $S$ containing an idempotent element, $C(S)$ is the minimum subsemigroup of $S$ with $E(S)$ as its set of idempotents.
Recall that if $S$ is regular then so is $C(S)$.

The concept of a biordered set comes into the center of interest outside regular semigroups when Easdown \cite{Easdown} 
proves that the axioms of a biordered set fully characterize the biordered sets of arbitrary semigroups.
In the proof he constructs, for any biordered set $E$, an idempotent generated semigroup with biordered set $E$.
Namely, he considers the semigroup $\IG(E)$ given by the presentation 
$\langle E\,|\, R\rangle$ 
where $R$ consists of all equalities $ef=g$ of words that are valid in $E$ for the partial operation,
and he shows that the natural mapping sending each element of $E$ to the element of $\IG(E)$ it represents is an 
isomorphism from $E$ to the biordered set of $\IG(E)$, which allows one to identify them.
The semigroup $\IG(E)$ is called the {\em free idempotent generated semigroup on $E$} since it has the universal property
that, whenever $S$ is a semigroup with biordered set $E$ such that $S=C(S)$, the identity mapping $E\to E$ 
can be extended to a homomorpism $\IG(E)\to S$. Clearly, this homomorphism is uniquely determined and surjective.

The existence of the regular analogue of $\IG(E)$, that is, of the free regular idempotent generated semigroup $\RIG(E)$
on a regular biordered set $E$, is proved by Nambooripad \cite{Nambooripad-79} (see also \cite[Section IX.1]{Grillet-95}).
Note that $\RIG(E)$ can be also obtained by presentation if the equalities $egf=ef$ are added to $R$ for every $e,f\in E$ and $g\in \mcS(e,f)$.

The question naturally arises how to recognize the biordered sets of important subclasses of semigroups.
A biordered set $E$ is the biordered set of a locally inverse semigroup if and only if $\mcS(e,f)$ is a singleton for every
$e,f\in E$.
In this case, the algebra $(E,\wedge)$ where $e\wedge f$ is the unique element of $\mcS(f,e)$ is called a 
{\em pseudosemilattice}.
Nambooripad \cite{Nambooripad-81} proves that an algebra $(E,\wedge)$ is a pseudosemilattice if and only if it satisfies certain identities.
The first description of the free pseudosemilattice on a set is provided by Meakin \cite{Meakin-83}, and certain subvarieties
of the variety of pseudosemilattices are investigated by Auinger and Oliveira \cite{Auin-Oliv}.
The latter paper gives an alternative model for the free pseudosemilattice on a set, and refers to several papers published in the meantime which also contain models of different kinds. 

It is proved by Clifford \cite{Clifford-76} that a regular biordered set comes from a completely regular semigroup 
(equivalently, from an $E$-solid semigroup) if and only if the relations $\mcLx$ and $\mcRx$ commute.
An alternative way, similar to pseudosemilattices, is found by Broeksteeg \cite{Broeksteeg} to describe these biordered sets where he associates an algebra
$(E,\ast)$ to the biordered set $E$ of every completely regular semigroup $S$ by defining $e\ast f$ to be the idempotent 
element of the $\mcHx$-class of $ef$ in $S$, and he axiomatizes these algebras by identities.
The biordered sets of bands are characterized in \cite{Nambooripad-79}, 
and a more intrinsic description of them is due to Easdown \cite{Easdown-84a}.
The biordered sets of eventually regular semigroups and those of certain special eventually regular semigroups (e.g., group-bound, periodic and finite semigroups) are described by Easdown \cite{Easdown-84b}.

Considerable attention has been devoted to free idempotent generated semigroups, and a number of deep results have been achieved.
If $E$ is a (regular) biordered set and $S$ is a (regular) idempotent generated semigroup with biordered set $E$
then the $\mcDx$-class of $S$ and that of $\IG(E)$ ($\RIG(E)$) containing $e\in E$ are similar to each other 
(e.g., the numbers of $\mcRx$-classes coincide, the same holds for $\mcLx$-classes).
Certainly, $\IG(E)$ might have non-regular $\mcDx$-classes. Moreover, the maximal subgroup of $S$ with identity element $e$
is a homomorphic image of the maximal subgroup of $\IG(E)$ with identity element $e$.
Therefore it is not surprising that the investigations have focused on the maximal subgroups of free (regular) idempotent generated semigroups.
Nambooripad and Pastijn prove in \cite{Namb-Past} that the maximal subgroups of $\RIG(E)$ are free for every 
pseudosemilattice $E$.
A partial generalization of this result for the non-regular case is published by McElwee in \cite{McElwee} 
where the maximal subgroups of $\IG(E)$ are shown to be free if the principal ideals of the biordered set $E$ are singletons. 
By these examples, it seemed to be plausible that this might be the case in general.
However, Brittenham, Margolis and Meakin disprove this expectation in \cite{Bri-Mar-Mea}.
First, they notice that, for every regular biordered set $E$, the natural homomorphism from $\IG(E)$ to $\RIG(E)$, when restricted to the set of regular elements of $\IG(E)$, is a bijection (although not an isomorphism since regular elements do not form a subsemigroup in $\IG(E)$ in general), and consequently, the maximal subgroups of $\IG(E)$ and of $\RIG(E)$ containing $e\in E$ are isomorphic.
The main result of the paper constructs a finite regular biordered set $E$ such that $\RIG(E)$ has the free abelian group of rank $2$ as a maximal subgroup.
The proof combines Nambooripad's theory \cite{Nambooripad-79} with topological methods.
Due to Gray and Ru\v{s}kuc \cite{Gray-Rusk}, it turns out that the case is just the opposite of the former expectations, 
since every group $G$ is a maximal subgroup of $\IG(E)$ for some biordered set $E$.
What is more, if $G$ is finitely presented then $E$ can be chosen to be a biordered set of a finite semigroup.
Dolinka and Ru\v{s}kuc \cite{Doli-Rusk} strengthen these results by proving that $E$ can be chosen in them to be a biordered set of a band.
On the other hand, Gould and Yang \cite{Gould-Yang} find a natural proof for the fact that any (finite) group is a maximal subgroup of some free idempotent generated semigroup over a (finite) biordered set which involves much less machinery
than the proof in \cite{Gray-Rusk}.

In a recent paper Dolinka, Gray and Ru\v{s}kuc \cite{Dol-Gra-Rus} initiate the study of algorithmic problems related to $\IG(E)$ provided $E$ is finite.
They present an algorithm which computes, for any $e\in E$, a finite presentation for the maximal subgroup of $\IG(E)$ contaning $e$, and they show that if all these finitely presented groups have decidable word problem then, for any pair of words in $E^*$, it is decidable whether they represent regular elements of $\IG(E)$, and if so, whether they represent the same element. 
However, they also prove that there exists a band such that $\IG(E)$ with $E$ being the biordered set of this band
has undecidable word problem.

\subsection{Generalization} \label{gen-fund}

Nambooripad's approach (Section \ref{Namboo-fund}) is partly generalized to the class of all semigroups by Easdown \cite{Easdown-88} in the following way.
Edwards \cite{Edwards-83, Edwards-85} introduced an idempotent separating congruence $\muc$ on any 
semigroup $S$ such that the $\muc$ relation of $S/\muc$ is the equality.
The relation $\muc$ need not be the maximum idempotent separating congruence on $S$ but it is so provided $S$ is eventually regular.
An important property of $\muc$ presented by Easdown and Hall \cite{Eas-Hall} is that it is biorder preserving, 
i.e., for every
$e,f\in E(S)$, we have $e\oml f$ if and only if $e\muc\oml f\muc$, and the dual property also holds.
Define a {\em fundamental semigroup} to be a semigroup whose only idempotent separating congruence is the equality relation,
and notice that this notion of a fundamental semigroup is more restrictive than that in \cite[Section III.4]{Grillet-95}. 
Easdown \cite{Easdown-88} proves that every semigroup where the relation $\muc$ is the equality is necessarily fundamental. Combining this fact with the previous ones, the main result follows: 
each semigroup is a biorder preserving extension by a fundamental semigroup.
Note that if $S$ is a biorder preserving extension by $T$ then $E(S)$ is, in general, only a biordered subset of $E(T)$.
It is also shown in \cite{Easdown-88} that the condition `biorder preserving extension' in the main result 
cannot be replaced by `extension which preserves the biordered set'.

The fundamental semigroups with a given biordered set are described by Easdown, Jordan and Roberts \cite{Eas-Jor-Rob} for the class of semigroups generated by regular elements.
For a biordered set $E$, a fundamental semigroup $T_E$ with biordered set $E$ is introduced which is generated by regular elements, and has the properties that the so-called symmetric subsemigroups (among them $T_E$ itself) are fundamental,
and that each semigroup $S$ generated by regular elements and having $E$ as its biordered set admits a biorder preserving representation $S\to T_E$ whose range is a symmetric subsemigroup of $T_E$.

\section{The structure of regular semigroups via groupoids and via cross-connections} \label{sect:grd-crc}

We continue our survey with one of the most substantial results in the theory of regular semigroups, also achieved in the 1970's: Nambooripad's description of the structure of regular semigroups via ordered groupoids over biordered sets.
Section \ref{Namboo-rig} is devoted to this highly non-trivial generalization of Schein's result \cite{Schein}
(see also \cite[Propositions 4.1.1 and 4.1.7]{Lawson}),
The concept taking over the role of an inductive groupoid is also called `inductive groupoid' by Nambooripad.
In order to avoid confusion with the notion mentioned earlier in the context of inverse semigroups,
we will use the expression `regular inductive groupoid' for Nambooripad's concept. 
The subject of Section \ref{Namboo-cc} is an alternative description of the structure of regular semigroups due to Nambooripad which generalizes Grillet's approach (Section \ref{Grillet}) for all regular semigroups.
Additionally, we discuss a recent work by Muhammed and Volkov \cite{Muha-Volk-19a,Muha-Volk-19b} on the relationship 
between Nambooripad's two approaches.
We close this section by outlining a generalization of Nambooripad's description of regular semigroups via 
ordered groupoids over biordered sets for a class of semigroups with distinguished sets of idempotents.

The idea behind Schein's result is fairly natural and transparent;
see \cite[Chapter 1]{Lawson} where the historical background is also cleared up.
There are two natural ways to define composition of partial bijections of a given set $X$: 
either we compose them as partial transformations, in which case product is defined for any pair $\alpha,\beta$ of partial 
bijections on $X$, or we compose them as usual mappings, 
in which case product is defined if and only if the range of $\alpha$ coincides with the domain of $\beta$.
To distinguish the two kinds of product, the latter is called reduced product. 
Given a collection $F$ of partial bijections on $X$ such that the inverse of $\alpha$ belongs to $F$ for every $\alpha\in F$, 
it is easy to see that $F$ is closed under composition of partial bijections if and only if 
$F$ is closed under usual composition of mappings, and
$F$ is closed under restriction to all intersections of possible domains (ranges).   
The notion of an inverse semigroup and that of an inductive groupoid describe, up to isomorphism, these sets $F$ in the first and in the second sense, respectively.
More formally, Schein's result
provides an alternative approach to inverse semigroups formulated in the language of groupoids.
Namely, the natural partial order and the restricted product of an inverse semigroup $S$ lead to an inductive groupoid
with the semilattice $E(S)$ as its partially ordered set of objects, and conversely, the composition of an inductive groupoid can be naturally extended to an everywhere defined pseudoproduct, thus forming an inverse semigroup from the inductive groupoid. 
Moreover, these two constructions are inverses of each other.
When generalizing these results for regular semigroups, Nambooripad also considered homomorphisms between regular semigroups, and formulated his result in the language of isomorphisms of categories.
The second part of the Ehresmann--Schein--Nambooripad theorem \cite[Theorem 4.1.8]{Lawson} on the equivalence of the category of inverse semigroups and usual homomorphisms and of the category of inductive groupoids and inductive functors is the inverse semigroup version of his result.
In the sequel, when referring to the Ehresmann--Schein--Nambooripad theorem, we always mean this second part.

\subsection{Regular inductive groupoids} \label{Namboo-rig}

Originally, Nambooripad described the structure of all regular semigroups in \cite{Nambooripad-73} (see also \cite{Nambooripad-75b}) in terms of restricted product and mappings between certain $\mcRx$-classes and $\mcLx$-classes.
An alternative description was published by Meakin \cite{Meakin-78} where he used the warp of a regular semigroup rather than its biordered set.
In \cite{Nambooripad-79} Nambooripad reformulates his approach by making use of ordered groupoids over biordered sets.
He introduces the notion of a regular inductive groupoid and that of a regular inductive functor, and
he proves that the category of regular semigroups and usual homomorphisms and the category of regular inductive groupoids and regular inductive functors are equivalent.
Similarly to the fundamental case, Nambooripad's notion of a regular inductive groupoid is an exact generalization of the notion of an inductive groupoid.

For any regular semigroup $S$, an ordered groupoid $\mcG(S)$ is defined whose structure of objects is the biordered set
$E(S)$, and whose morphisms are the pairs $(s,s')$ where $s\in S$ and $s'\in V(s)$, and we have $\dd(s,s')=ss'$ and 
$\rr(s,s')=s's$ (with the respective identity morphisms being $(ss',ss')$ and $(s's,s's)$, respectively).
Composition of morphisms $(s,s'),(t,t')$ is defined to be $(st,t's')$ provided $\rr(s,s')=\dd(t,t')$, and the inverse of
$(s,s')$ is $(s',s)$.
The partial order of $\mcG(S)$ is given by the rule
\[(s,s')\le (t,t')\quad \hbox{if and only if}\quad s=(ss')t,\ s'=t'(ss')\ \hbox{and}\ ss'\om tt'.\]

If $S$ is an inverse semigroup then $\mcG(S)$ is isomorphic to the inductive groupoid associated to $S$ by the Ehresmann--Schein--Nambooripad theorem, and the ordered groupoids $\mcG(S)$ of inverse semigroups $S$ can be characterized by the property that their sets of objects are semilattices.
With regular semigroups, the case is much more complicated.
The family of ordered groupoids with biordered sets as sets of objects is not sufficient to describe the structure of regular semigroups (\cite{Nambooripad-X18}).
Notice that if $S$ is an inverse semigroup then the morphisms of $\mcG(S)$ are in one-to-one correspondence with the elements of $S$, but this is by no means the case with regular semigroups.
If $S$ is any regular semigroup then, in order to recapture $S$ from $\mcG(S)$, one needs to understand and describe how these morphisms might relate to each other in $S$.
The additional structure introduced by Nambooripad to obtain the notion of a regular inductive groupoid is expressed in terms of an ordered groupoid constructed from sequences of elements of a regular biordered set.

Let $E$ be a regular biordered set. An \emph{$E$-path} is a sequence $(e_1,\ldots,e_n)$ of elements of $E$ where
$e_i\,(\mcLx\cup\mcRx)\,e_{i+1}$ for $i=1,\ldots,n-1$. 
An element $e_i$ is inesseantial if $e_{i-1}\mcL e_i\mcL e_{i+1}$ or $e_{i-1}\mcR e_i\mcR e_{i+1}$. 
The relation on the set of $E$-paths obtained by adding or removing inessential elements is an equivalence relation, and equivalence classes of $E$-paths are called \emph{$E$-chains}. 
Note that in each $E$-chain, there is a unique $E$-path without inessential elements, and it is used to represent the $E$-chain.
The $E$-chains form the morphisms of the ordered groupoid $\mcC(E)$ whose set of objects is $E$, and an $E$-chain 
$(e_1,\ldots,e_n)$ is a morphism from $e_1$ to $e_n$. 
Product of composable $E$-chains $(e_1,\ldots,e_n)$ and $(f_1,\ldots,f_r)$ is defined to be the $E$-chain of the $E$-path $(e_1, \ldots, e_n=f_1, \ldots, f_r)$, 
and the partial order of $E$-chains is defined as follows: 
$(e_1, \ldots, e_n) \le (f_1, \ldots, f_r)$ if and only if $e_1 \om f_1$ and the $E$-chain $(e_1, \ldots, e_n)$ contains the $E$-path $(e_1=f_1e_1f_1, f_2e_1f_2, f_3f_2e_1f_2f_3, \ldots, f_r\cdots f_2e_1f_2\cdots f_r)$.

Nambooripad introduced regular inductive groupoids in the following way.
Given a biordered set $E$ and an ordered groupoid $\mcG$, an \emph{evalution functor from $\mcC(E)$ to $\mcG$} is a 
functor $\eps\colon \mcC(E)\to \mcG$ whose object mapping is an order isomorphism.
The pair $(\mcG,\eps)$ is called a \emph{regular inductive groupoid} if it satisfies two axioms and their duals.
The main result on the general structure of a regular semigroup is the following.

\begin{theorem}[Nambooripad \cite{Nambooripad-79}] \label{Namb-reg}
The category of regular semigroups with usual homomorphisms as morphisms and the category of regular inductive groupoids with  regular inductive functors as morphisms are equivalent.
\end{theorem}

To give some insight into the proof, we mention that the regular inductive groupoid assigned to a regular semigroup $S$ with biordered set $E$ is $(\mcG(S),\eps_S)$ where $\mcG(S)$ is defined above, and $\eps_S$ is given by the rule 
$(e_1, \ldots, e_n)\mapsto (e_1 \cdots e_n, e_n \cdots e_1)$.
In the reverse direction, if $(\mcG,\eps)$ is a regular inductive groupoid then an equivalence relation $\equiv$ is defined
on the set of morphisms of $\mcG$ as follows:
\[x\equiv y\quad \hbox{if and only if}\quad \dd(x)\mcR \dd(y),\ \rr(x)\mcL \rr(y)\ \,\hbox{and}\,\ 
x\eps(\rr(x),\rr(y))=\eps(\dd(x),\dd(y))y.\]
On the set of $\equiv$-classes, an appropriate multiplication can be introduced such that a regular semigroup is obtained.
This regular semigroup is assigned to $(\mcG,\eps)$.

It is important to note that, restricting the constructions to inverse semigroups, Theorem \ref{Namb-reg} specializes to the
Ehresmann--Schein--Nambooripad theorem.

\subsection{Cross-connections} \label{Namboo-cc}

In a subsequent work \cite{Nambooripad-89} (see also \cite{Nambooripad-94}), Nambooripad extends Grillet's cross-connection approach from fundamental regular semigroups to arbitrary regular semigroups, again making use of categories, but this time the tools needed are deeper and more technical than before.
This might be one of the reasons that these results have aroused much less interest than his earlier work.
Another reason might be that 
\cite{Nambooripad-89} and \cite{Nambooripad-94} appeared as local publications and have remained hidden for the research community. 
Roughly speaking, in his construction, Nambooripad replaces Grillet's regular partially ordered sets by categories.
More precisely, in order to describe the structure of principal one-sided ideals of an arbitrary regular semigroup in an abstract way, he introduces the notion of a normal category, and, to model the interrelations between the principal left and principal right ideals, he generalizes Grillet's notion of cross-connection to normal categories so that cross-connections form a category (with respect to appropriate morphisms). 
He obtains the following result.

\begin{theorem}[Nambooripad \cite{Nambooripad-94}] \label{Namb-cross}
The category of regular semigroups with usual homomorphisms as morphisms and the category of cross-connections with appropriate functors as morphisms are equivalent.
\end{theorem}

Notice that, by transitivity, this theorem combined with Theorem \ref{Namb-reg} implies that 
the the category of cross-connections and the category of regular inductive groupoids are also equivalent.

Motivated by Theorems \ref{Namb-reg} and \ref{Namb-cross}, Muhammed and Volkov have investigated the interrelations between cross-connections and regular inductive groupoids.
In \cite{Muha-Volk-19a} they construct the regular inductive groupoid of a regular semigroup directly from
the cross-connection representation of the semigroup, and vice versa, by analyzing the rather complicated relationship
between the idempotent structure and ideal structure of an arbitrary regular semigroup.
It should be emphasized that the equivalence of categories established is not the composition of the equivalences given in the proofs of Theorems \ref{Namb-reg} and \ref{Namb-cross}.
In \cite{Muha-Volk-19b} the same authors go further by providing an equivalence between the category of regular inductive groupoids and the category of cross-connections in such a way that they avoid using semigroups and restrict themselves to a purely categorical framework.
It is worth mentioning that cross-connections seem to encode much more information than regular inductive groupoids.
Given a cross-connection, the regular inductive groupoid assigned to it by the equivalence presented in the paper can be found in a fairly straightforward and transparent way `inside' the cross-connection.
However, in the reverse direction, the construction of the cross-connection corresponding to a regular inductive groupoid is rather complicated: first, several auxiliary categories are extracted from the regular inductive groupoid, and then, these categories are combined to build the ingredients of the cross-connection.

\subsection{Generalization}

The most general class of semigroups for which Nambooripad's regular inductive groupoid approach is generalized is a class of semigroups with distinguished subsets of idempotents.

Let $S$ be a semigroup and $U$ a non-empty subset of $E(S)$.
The pair $(S,U)$ is said to be a {\em semigroup $S$ with distinguished subset of idempotents $U$}.
Consider the following relation $\tmcLUx$ on $S$, and its dual $\tmcRUx$, called generalized Green's relations: 
for every $a,b\in S$, let $a\tmcLU b$ if and only if the equalities $ae=af$ and $be=bf$ hold for the same pairs 
of idempotents $e,f\in U$.
These relations are known to be equivalence relations such that $\mcLx\subseteq \tmcLUx$ and $\mcRx\subseteq \tmcRUx$.
Note that $\tmcLUx$ ($\tmcRUx$) need not be a right (left) congruence.
Given semigroups $S$ and $T$ with distinguished subsets of idempotents $U$ and $V$, respectively, a homomorphism
$\phi\colon S\to T$ is said to be {\em admissible} if $U\phi\subseteq V$, and the following property and its dual are valid:
for every $a,b\in S$, the relation $a\tmcLU b$ implies $a\phi\tmcLV b\phi$.
For a general introduction to the topic of semigroups with distinguished subsets of idempotents, see Gould \cite{Gould-restr}.

A semigroup $S$ with a distinguished subset of idempotents $U$ is called in Wang \cite{Wang} {\em weakly regular} (or briefly and somewhat unprecisely, $S$ is called weakly $U$-regular) if the following three conditions are satisfied:
every $\tmcLUx$-class and every $\tmcRUx$-class contains an element of $U$,
$\tmcLUx$ is a right congruence and $\tmcRUx$ is a left congruence on $S$, and
the subsemigroup $\langle U\rangle$ generated by $U$ is a regular subsemigroup in $S$ with $E(\langle U\rangle)=U$.
Clearly, the definition implies that $U$ forms a regular biordered set.
This allows us to introduce the more informative name {\em weakly regular semigroup with a distinguished biordered subset} for these structures.
Notice that each regular semigroup $S$ is a weakly regular semigroup with distinguished biordered subset $E(S)$.

Theorem \ref{Namb-reg} is extended to weakly regular semigroups with distinguished biordered subsets by Wang \cite{Wang} as follows. 
The notion of a weakly regular category over a regular biordered set is introduced, and it is proved
that the category of weakly regular semigroups with distinguished biordered subsets together with admissible homomorphisms as morphisms and the category of weakly regular categories over regular biordered sets with appropriate functors as morphisms are equivalent.

The result by Easdown, Jordan and Roberts \cite{Eas-Jor-Rob} outlined in Section \ref{gen-fund} suggests the potentiality of generalizing Theorem \ref{Namb-reg} for the class of all semigroups generated by their regular elements, and perhaps even further, for a class of semigroups with distinguished biordered sets which contains both these semigroups and the weakly regular semigroups.

Before closing this section, we call the attention to \cite[Section 6]{Muha-Volk-19a} where
the authors formulate a bunch of topics for further research, mainly ones related to cross-connections.

\section{Structure theorems for classes of reguar semigroups motivated by the McAlister--Lawson theory} \label{sect:McA-Law}

While the Munn representation and the Ehresmann--Schein--Nambooripad theorem have exact generalizations for the whole class of regular semigroups, this is far not the case with the third substantial structure theory for inverse semigroups due to McAlister and Lawson (\cite[Sections 2.2, 7.1 and 7.2]{Lawson}).
Almost all elements of this theory have been generalized for orthodox and for locally inverse semigroups.
We summarize these results in Sections \ref{O-sd} and \ref{LI-Pp}, respectively.
The central point of the McAlister--Lawson theory is the class of $E$-unitary inverse semigroups, or alternatively, the class of extensions of semilattices by groups.
This has naturally led to the idea of applying semidirect products of semilattices by groups to describe their structure, and motivated the study of embeddability of an extension by a group and by an inverse semigroup, in general, in a semidirect or a $\lambda$-semidirect product.
For example, each idempotent separating extension by an inverse semigroup $T$, which plays essential role in Munn's approach for constructing all inverse semigroups (see the introduction of Section \ref{sect:fund}), is known to be embeddable in a $\lambda$-semidirect product of a Clifford semigroup by $T$ (\cite[Theorem 5.3.5 and Proposition 5.3.6]{Lawson}).
Section \ref{ext-sd} contains embedding theorems of this kind.
Section \ref{McL-gen} is mainly about a generalization of McAlister's covering theorem for arbitrary semigroups.

Motivated by Scheiblich's model for free inverse semigroups (\cite[Section 6.5]{Lawson}), 
McAlister proved his well-known theorems on $E$-unitary inverse semigroups 
(\cite[Theorems 2.2.4, 7.2.15]{Lawson}). 
By McAlister's covering theorem, each (finite) inverse semigroup has a (finite) $E$-unitary cover,
and his $P$-theorem describes the structure of a (finite) $E$-unitary inverse semigroup 
by means of a (finite) partially ordered set, a (finite) semilattice and a (finite) group with a construction reminicent to a semidirect product of a semilattice by a group.
O'Carroll applied this structure theorem to show that each (finite) $E$-unitary inverse semigroup is embeddable in a (finite) semidirect product of a semilattice by a group (\cite[Theorem 7.1.5]{Lawson}).
Note that this property characterizes (finite) $E$-unitary inverse semigroups since an inverse subsemigroup of a semidirect product of a semilattice by a group is easily seen to be $E$-unitary.
Combining the covering and the embedding theorems, one immediately obtains that 
every (finite) inverse semigroup divides a (finite) semidirect product of a semilattice by a group 
(\cite[Theorem 7.1.6]{Lawson}).
This division theorem constructs inverse semigroups from semilattices and groups in a way which significantly differs from that implied by the Munn representation and by the description of fundamental inverse semigroups.

An elegant proof of McAlister's covering theorem is based on the result that each (finite) inverse semigroup can be embedded into a (finite) factorizable inverse monoid (\cite[Theorem 2.2.3]{Lawson}).
It is fairly easy to see that every factorizable inverse monoid is an idempotent separating homomorphic image of a semidirect product of its semilattice (which is a monoid) by its group of units (\cite[Corollary 7.1.11]{Lawson}).
Investigating the homomorphic images of semidirect products of semilattices by groups in general, McAlister and Lawson developed an alternative way for producing inverse semigoups from semidirect products of semilattices by groups.
They introduced the notion of an almost factorizable inverse semigroup, characterized them as the inverse semigroups obtained from factorizable inverse monoids by removing their groups of units
(\cite[Proposition 7.1.12 and Theorem 7.1.13]{Lawson}), and proved that
an inverse semigroup is almost factorizable if and only if it is a homomorphic image (equivalently, an idempotent separating homomorphic image) of a semidirect product of a semilattice by a group
(\cite[Theorem 7.1.10 and its proof]{Lawson}).
Finally, it is worth mentioning that the intersection of the class of $E$-unitary and the class of almost factorizable inverse semigroups is just the class of semidirect products of semilattices by groups
(\cite[Theorem 7.1.8]{Lawson}).

\subsection{Orthodox semigroups and semidirect products of bands by groups} \label{O-sd}

A regular semigroup $S$ is said to be {\em $C_\infty$-unitary} if $C_\infty(S)$ is a unitary subset in $S$, or, equivalently, $S$ has a group congruence with kernel $C_\infty(S)$.
The latter property implies that this group congruence is necessarily the minimum group congruence $\sigma$ on $S$.
If $S$ is orthodox then $C_\infty(S)=E(S)$, and so this definition yields the notion of an {\em $E$-unitary orthodox semigroup} which obviously reduces to the usual notion of $E$-unitariness in the inverse case.
In the literature, $E$-unitary orthodox semigroups are frequently called $E$-unitary regular semigroups because,
replacing $C_\infty(S)$ by $E(S)$ in the definition, the semigroups obtained are just the $E$-unitary orthodox semigroups.
In this survey we avoid using the term `$E$-unitary regular semigroup'.
A homomorphism $\phi\colon S\to T$ between regular semigroups is said to be {\em $C_\infty$-separating} if the restriction of $\phi$ to $C_\infty(S)$ is injective. 

In this context, McAlister's covering theorem extends to all regular semigroups, 
an important intermediate step being the case of orthodox semigroups 
(\cite[Theorem IX.4.2]{Grillet-95}).

\begin{theorem}[Trotter~\cite{Trotter-95}] \label{cov-Cinf}
Each (finite) regular semigroup is a $C_\infty$-separating homomorphic image of a (finite) $C_\infty$-unitary regular semigroup.
\end{theorem}

It is clear by definition that a regular semigroup $S$ is $C_\infty$-unitary if and only if it is an extension of 
$C_\infty(S)$ by the group $S/\sigma$.
Therefore a generalization of O'Carroll's embedding theorem amounts to studying which classes $\mathcal{K}$ of regular semigroups have the property that if $S$ is a $C_\infty$-unitary regular semigroup with $C_\infty(S)\in \mathcal{K}$ then
$S$ is embeddable in a semidirect product of a member of $\mathcal{K}$ by a group.

O'Carroll's embedding theorem has been generalized by the author for the class of $E$-unitary orthodox semigroups whose bands are regular, i.e., whose bands satisfy the identity $axya=axaya$
(\cite[Theorem IX.5.7]{Grillet-95}).
Note that it is not emphasized in the statement that the finitary version also holds but it follows easily from the construction in the proof (see the paper cited), and from the well-known fact that finitely generated bands are finite.

The restriction on the bands of $E$-unitary orthodox semigroups cannot be removed from \cite[Theorem IX.5.7]{Grillet-95}.
A finite $E$-unitary orthodox semigroup is constructed by Billhardt \cite{Billhardt-98} which fails to be embeddable in a semidirect product of a band by a group, and whose band is left semiregular, i.e., belongs to a variety being `close' to the variety of regular bands in the lattice of band varieties.
However, the division theorem extends to the class of all orthodox semigroups (\cite[Theorem IX.5.8]{Grillet-95}), although its finitary version is left open.

Homomorphic images of semidirect products of bands by groups are studied by Hartmann \cite{Hartm-07}.
For orthodox monoids, the notion of factorizability and, for orthodox semigroups, the notion of almost factorizability can be generalized in a way that they relate to each other in the same way as in the inverse case.
The connection between factorizable orthodox monoids and semidirect products of band monoids by groups
is also the same as in the inverse case. 
However, the case of almost factorizable orthodox semigroups and semidirect products of bands by groups turns out to be more complicated, and only the following weaker result holds.

\begin{theorem}[Hartmann \cite{Hartm-07}]
An orthodox semigroup is almost factorizable if and only if it is an idempotent separating homomorphic image of a semidirect product of a band by a group.
\end{theorem}

The class of almost factorizable orthodox semigroups is shown to be a proper subclass of the class of all homomorphic images of semidirect products of bands by groups.
It is also established that a generalized inverse semigroup $S$ is a homomorphic image of a semidirect product of a band (or, equivalently, of a normal band) by a group if and only if the maximum inverse semigroup quotient of $S$ is almost factorizable.
In general, almost factorizability of the maximum inverse semigroup quotient is necessary but not sufficient for an orthodox semigroup to be a homomorphic image of a semidirect product of a band by a group.

After having these facts at hand, it is not surprising that the class of semidirect products of bands by groups is a proper subclass of the intersection of the class of $E$-unitary orthodox semigroups and the class of almost factorizable orthodox semigroups. 
To demnostrate this, Hartmann and the author \cite{Hart-Sz} give a finite example which is very close to being inverse, namely, it is a finite generalized inverse semigroup whose band is left normal.
Furthermore, they characterize the $E$-unitary and almost factorizable orthodox semigroups that are (isomorphic to) semidirect products of bands by groups, and by making use of a construction generalizing semidirect products, they describe the structure of $E$-unitary and almost factorizable orthodox semigroups by means of bands and groups.

\subsection{Locally inverse semigroups and Pastijn products of normal bands by completely simple semigroups} \label{LI-Pp}

Let $S$ and $T$ be regular semigroups.
A subsemigroup of $S$ of the form $eSe\ (e\in E(S))$ is called a {\em local subsemigroup} of $S$.
A homomorphism $\phi$ from $S$ onto $T$ is said to be a {\em local isomorphism} if $\phi$ is one-to-one on every local subsemigroup of $S$. 
A {\em regular Rees matrix semigroup over a regular semigroup $S$} is defined to be the subsemigroup of all regular elements of a Rees matrix semigroup $M(S;I,\Lambda;P)$ over $S$, and is denoted $R(S;I,\Lambda;P)$.
A widely known structure theorem, due to McAlister (\cite[Section IX.3]{Grillet-95}) characterizes locally inverse semigroups in terms of regular Rees matrix semigroups, and an analogous result has been published by the same author for further classes.

\begin{theorem}[McAlister \cite{McAlister-84}]
A regular semigroup is locally inverse (locally orthodox, locally $E$-solid) if and only if it is a locally isomorphic image of a regular Rees matrix semigroup over an inverse (orthodox, $E$-solid) semigroup.
\end{theorem}

Applying the argument due to McAlister \cite{McAlister-83} in the locally inverse case, this theorem combined with 
Theorem \ref{cov-Cinf} implies the following covering theorem.

\begin{theorem}
For every locally inverse (locally orthodox, locally $E$-solid) semigroup $S$, there exists a regular Rees matrix semigroup 
$R$ over an $E$-unitary inverse ($E$-unitary orthodox, $C_\infty$-unitary $E$-solid) semigroup such that $S$ is a homomorphic image of $R$ over completely simple semigroups.
\end{theorem}

A locally inverse semigroup where the set of idempotents is a disjoint union of semilattices is called {\em straight}
(Pastijn--Petrich \cite{Pastijn-Petrich-84}).
We say that a locally inverse semigroup $S$ is {\em right straight} if it is a rectangular band $I\times \Lambda$ of its subsemigroups $S_{i\lambda}\ ((i,\lambda)\in I\times \Lambda)$ such that 
$E_i=\bigcup_{\lambda\in \Lambda} E(S_{i\lambda})$ is a right normal band for every $i\in I$, 
and a left straight locally inverse semigroups is defined dually.
One can see that a locally inverse semigroup is straight if and only if it is right and left straight.

The notion of a {\em weakly $E$-unitary locally inverse semigroup $S$} is introduced by Veeramony \cite{Veeramony} by requiring the same property of the natural partial order for a locally inverse semigroups which defines $E$-unitary inverse semigroups, i.e., that
$e\le a$ implies $a\in E(S)$ for every $e\in E(S)$ and $a\in S$.
For example, it is easy to see that a regular Rees matrix semigroup over an $E$-unitary inverse semigroup is a straight weakly $E$-unitary locally inverse semigroup.
It is proved by Ka\softd{}ourek \cite{Kadourek} that a locally inverse semigroup is weakly $E$-unitary if and only if its minimum completely simple congruence is idempotent pure.
Alternatively, the weakly $E$-unitary locally inverse semigroups are just the locally inverse semigroups which are extensions by completely simple semigroups over normal bands.

A construction producing a weakly $E$-unitary locally inverse semigroup from a normal band and a completely simple semigroup is due to Pastijn \cite{Pastijn-82}. 
It is also applied in \cite{Kadourek} and it is given the name Pastijn product.
Let $N$ be a normal band and let $T=M[G;I,\Lambda;P]$ be a completely simple semigroup represented as a Rees matrix semigroup where the group $G$ acts on $N$ by automorphisms.
Define a multiplication on $K\times T$ by the rule
\[\left(a,(i,g,\lambda)\right)\left(b,(j,h,\mu)\right)=
\left(a\cdot\acth{gp_{\lambda j}}b,(i,gp_{\lambda
j}h,\mu)\right).\]
The weakly $E$-unitary locally inverse semigroup so obtained is called a {\em Pastijn product of $N$ by $T$} and is denoted $N\odot T$.
In particular, if $T=G$ then $N\odot G$ is the semidirect product of $N$ by $G$, so Pastijn products generalize semidirect products of semilattices by groups.  

The next theorem generalizes O'Carroll's embedding theorem for weakly $E$-unitary locally inverse semigroups.

\begin{theorem}[Billhardt and Szendrei \cite{Billhardt-Sz}]
Each (finite) weakly $E$-unitary locally inverse semigroup $S$ can be embedded into a (finite) Pastijn product of a normal band $B$ by a completely simple semigroup. 
In particular, if $S$ is straight then $B$ can be chosen to be a semilattice, and if $S$ is left (right)
straight then $B$ can be chosen to be a left (right) normal band.
\end{theorem}

Homomorphic images of Pastijn products are studied by the author in \cite{Sz-X11}.
A locally inverse semigroup $S$ is called {\em almost factorizable} if there exists a completely simple subsemigroup
$\mcU$ in the semigroup $\mcO(S)$ of all order ideals of $S$ such that $\bigcup E(\mcU)=E(S)$ and 
$\bigcup \mcU\supseteq S$.
The result obtained is an exact generalization of the characterization of almost factorizable inverse semigroups in terms of semidirect products of semilattices by groups.

\begin{theorem}[Szendrei \cite{Sz-X11}]
For any locally inverse semigroup $S$, the following statements are equivalent:
\begin{enumerate}
\item $S$ is almost factorizable,
\item $S$ is a homomorphic image of a Pastijn product of a normal band by a completely simple semigroup,
\item $S$ is a homomorphic image of a Pastijn product of a normal band by a completely simple semigroup over completely simple semigroups.
\end{enumerate}
\end{theorem}

\subsection{Extensions embeddable in ($\lambda$-)semidirect products} \label{ext-sd}

The theorem on the embeddability of an $E$-unitary orthodox semigroup whose band is regular (Section \ref{O-sd}) has been generalized for extensions of regular orthogroups by groups as follows.
By definition, a {\em regular orthogroup} is an orthodox completely regular semigroup whose band is regular.

\begin{theorem}[Szendrei \cite{Sz-95}]
Every extension of a regular orthogroup $K$ by a group $G$ is embeddable in a semidirect product of a regular orthogroup $K'$ by $G$ where $K'$ belongs to the variety of orthogroups generated by $K$.
\end{theorem}

An easy consequence of \cite[Theorem 5.3.5 and Proposition 5.3.6]{Lawson} is the following.

\begin{theorem}
For every inverse semigroup $S$ and idempotent separating congruence $\theta$, the extension $\sthe$ can be embedded into a $\lambda$-semidirect product of a Clifford semigroup $K$ by $S/\theta$ where $K$ belongs to the variety of Clifford semigroups generated by 
$\Ker \theta$.
\end{theorem}

Since the idempotent separating extensions by inverse semigroups are just the extensions by inverse semigroups over groups, and the kernel of a $\lambda$-semidirect product of a group by an inverse semigroup, as an extension, is a Clifford semigroup, the question arises whether idempotent separating extensions by inverse semigroups can be embedded into $\lambda$-semidirect products of groups by inverse semigroups.
The affirmative answer to this question strengthens the previous theorem.

\begin{theorem}[Billhardt and Szittyai \cite{Bill-Szitt}]
For every inverse semigroup $S$ and idempotent separating congruence $\theta$, the extension $\sthe$ can be embedded into a $\lambda$-semidirect product of a group $H$ by $S/\theta$ where $H$ belongs to the variety of groups generated by the $\theta$-classes containing idempotents.
\end{theorem}

Recently, a weakened version of this result (where $H$ is only required to be a group) has been generalized for certain extensions of completely simple semigroups by inverse semigroups so that a structure theorem for $E$-solid locally inverse semigroups is implied.
Recall (Section \ref{Prelim}) that a regular semigroup is $E$-solid if and only if
it is an extension by an inverse semigroup over completely simple semigroups.
A $\lambda$-semidirect product of a completely simple semigroup by an inverse semigroup is easily checked to be locally inverse.
Therefore one should restrict his attention to $E$-solid locally inverse semigroups when studying which extensions by inverse semigroups over completely simple semigroups are embeddable in $\lambda$-semidirect products of completely simple semigroups by inverse semigroups.

\begin{theorem}[D\'{e}k\'{a}ny, Szendrei and Szittyai \cite{Dek-Sz-Szit}]
Let $S$ be an $E$-solid locally inverse semigroup and $\theta$ an inverse semigroup congruence on $S$ over completely simple semigroups.
Then the extension $\sthe$ can be embedded into a $\lambda$-semidirect product of a completely simple semigroup by $S/\theta$.
\end{theorem}

Thus $E$-solid locally inverse semigroups can be built up from completely simple semigroups and inverse semigroups in terms of two fairly simple constructions: forming $\lambda$-semidirect product and taking regular subsemigroup.
Conversely, only $E$-solid locally inverse semigroups can be produced in this way.

\subsection{Generalization} \label{McL-gen}

Fountain, Pin and Weil \cite{Fountain-Pin-Weil} generalize Theorem \ref{cov-Cinf} for the class of arbitrary monoids (and most of the results can be adjusted also to arbitrary semigroups).
They develop a theory which describes extensions of monoids by groups in terms of groups acting on categories, and they 
apply it to find a sufficient condition for a monoid $M$ and its submonoid $T$ to possess a $T$-cover.
(The monoid version of Theorem \ref{cov-Cinf} is the case where $M$ is regular and $T=C_\infty(M)$.)
It is worth mentioning that this theory can be also applied to recover McAlister's $P$-theorem and most of the structure theorems published so far that describe members of wider classes as extensions by groups.

Notice that this theory says nothing about embeddability of extensions by groups in semidirect products, and the results on locally inverse semigroups presented in Section \ref{LI-Pp} are also outside the scope of this theory.
The author believes that the role of semidirect products of semilattices by groups is fundamental in the McAlister--Lawson approach to the structure of inverse semigroups, and so a generalization of this approach to wider classes of regular semigroups should contain an appropriate ingredient.
The results of Section \ref{LI-Pp} show that a generalization for a class containing locally inverse semigroups cannot be restricted so that the role of $E$-unitary inverse semigroups is taken over by semigroups which are extensions by group.
It seems to be a challenging research problem how to extend the McAlister--Lawson approach to the class of all regular semigroups, or at least, to a subclass containing both orthodox and locally inverse semigroups.

\end{document}